\documentclass[12pt]{amsart}
\newtheorem{theorem}{Theorem}[section]

\theoremstyle{definition}
\newtheorem{definition}{Definition}[section]
\newtheorem{corollary}{Corollary}[section]
\newtheorem{example}{Example}[section]

\theoremstyle{remark}

\numberwithin{equation}{section}

\begin{document}
\title[Ricci solitons on Ricci pseudosymmetric $(LCS)_n$-manifolds]{Ricci solitons on Ricci pseudosymmetric $(LCS)_n$-manifolds}
\author[S. K. Hui, R. S. Lemence and D. Chakraborty]{Shyamal Kumar Hui, Richard S. Lemence and Debabrata Chakraborty}
\subjclass[2000]{53B30, 53C15, 53C25}
%%%%%%%%%%%%%%%%%%%%%%%%%%%%%%55
\keywords{ Ricci soliton, $(LCS)_n$-manifold, Ricci pseudosymmetric manifold, concircular curvature tensor,
projective curvature tensor, $W_3$-curvature tensor, conharmonic curvature tensor, conformal curvature tensor.}
%%%%%%%%%%%%%%%%%%%%%%
\begin{abstract}
The object of the present paper is to study some types of Ricci
pseudosymmetric $(LCS)_n$-manifolds whose metric is Ricci soliton.
We found the conditions when Ricci soliton on concircular Ricci
pseudosymmetric, projective Ricci pseudosymmetric, $W_{3}$-Ricci
pseudosymmetric, conharmonic Ricci pseudosymmetric, conformal Ricci
pseudosymmetric $(LCS)_n$-manifolds to be shrinking, steady and
expanding.  We also construct an example of concircular
Ricci pseudosymmetric $(LCS)_3$-manifold whose metric is Ricci
soliton.
\end{abstract}
%%%%%%%%%%%%%%%%%%%%%%%%
\maketitle
%%%%%%%%%%%%%%%%%%%%%%%%%%%%%%%%%%%%%%%%%%%%%%%%%%%%%%%%%%%%%%%%%%%%%%%%
\section{Introduction}
%%%%%%%%%%%%%%%%%%%%%
In 2003, Shaikh \cite{SHAIKH1} introduced the notion of Lorentzian
concircular structure manifolds (briefly, $(LCS)_n$-manifolds) with
an example, which generalizes the notion of LP-Sasakian manifolds
introduced by Matsumoto \cite{MAT} and also by Mihai and Rosca
\cite{MIHAI}. Then Shaikh and Baishya (\cite{SHAIKH3}, \cite{SHAIKH4})
investigated the applications of $(LCS)_n$-manifolds to the general
theory of relativity and cosmology. The $(LCS)_n$-manifolds are also
studied by Atceken et. al. (\cite{ATCEKEN3}, \cite{ATEHUI}, \cite{SHMA}),
Hui \cite{SKH}, Hui and Chakraborty \cite{DEBA1},
Narain and Yadav \cite{DNSY}, Prakasha \cite{DGP},
Shaikh and his co-authors (\cite{SHAIKH11}, \cite{AHMAD}, \cite{SABINA1}--\cite{SB}, \cite{HUI}, \cite{SMHUI}) and many others.
%%%%%%%%%%%%%%%%%%%%%%%%%%%%%%%%%

%%%%%%%%%%%%%%%%%%%%%%%%%%%%%%%%%%%%5
In 1982, Hamilton \cite{HAM1} introduced the notion of Ricci flow to find a canonical metric on a smooth manifold. Then Ricci flow
has become a powerful tool for the study of Riemannian manifolds, especially for those manifolds with positive curvature.
Perelman (\cite{PER1}, \cite{PER2}) used Ricci flow and its surgery to prove the Poincare conjecture.
The Ricci flow is an evolution equation for metrics on a Riemannian manifold defined as follows:
%%%%%%%%%%%%%%%%%%%%%%%%%%%
\begin{equation*}
\frac{\partial}{\partial t}g_{ij}(t) = - 2 R_{ij}.
\end{equation*}
%%%%%%%%%%%%%%%%%%%%%%%%%%%%%%%%%%%%%%%%%%%%
\indent \indent A Ricci soliton emerges as the limit of the
solutions of the Ricci flow. A solution to the Ricci flow is called
Ricci soliton if it moves only by a one parameter group of
diffeomorphisms and scaling. To be precise, a Ricci soliton on a
Riemannian manifold $(M,g)$ is a triple $(g,V,\lambda)$ satisfying
\cite{HAM2}
\begin{equation}
\label{eqn1.1}
\pounds_Vg + 2S + 2\lambda g = 0,
\end{equation}
%%%%%%%%%%%%%%%%%%%%%%%
where $S$ is the Ricci tensor, $\pounds_V$ is the Lie derivative
along the vector field $V$ on $M$ and $\lambda \in \mathbb{R}$. The
Ricci soliton is said to be shrinking, steady and expanding
according as $\lambda$ is negative, zero and positive respectively.
%%%%%%%%%%%%%%%%%%%%%%%%%

%%%%%%%%%%%%%%%%%%%%%%%%%%%%%5
During the last two decades, the geometry of Ricci solitons has been
the focus of attention of many mathematicians. In particular, it has
become more important after Perelman applied Ricci solitons to solve
the long standing Poincare conjecture posed in 1904. In
\cite{Sharma3} Sharma studied the Ricci solitons in contact
geometry. Thereafter Ricci solitons in contact metric manifolds have
been studied by various authors such as Bagewadi et. al (\cite{Ash},
\cite{ABI}, \cite{BAG}, \cite{ING}), Bejan and Crasmareanu
\cite{Cras}, Blaga \cite{Bloga}, Chen and
Deshmukh \cite{CD}, Deshmukh et. al \cite{Detal}, Hui et. al
(\cite{CHS}, \cite{DEBA1}-\cite{SKH6}, \cite{DEBA2}-\cite{HUC}), Nagaraja and Premalatha \cite{NAGA},
Tripathi \cite{Trip} and many others.\\
%%%%%%%%%%%%%%%%%%%%%%%%%%%%%%%%%%%%%%%%%%%%%%%%%%%%%%%%%%%%
\indent The notion of Ricci pseudosymmetric manifold was introduced
by Deszcz (\cite{dez3},\cite{DESZCZ1}). A geometrical interpretation
of Ricci pseudosymmetric manifolds in the Riemannian case is given
in \cite{JHSV}. A $(LCS)_n$-manifold $(M^n,g)$ is called Ricci
pseudosymmetric (\cite{dez3}, \cite{DESZCZ1}) if the tensor $R\cdot
S$ and the Tachibana tensor $Q(g,S)$ are linearly dependent, where
%%%%%%%%%%%%%%%%%%%%%
\begin{equation}
\label{eqn1.2}
(R(X,Y)\cdot S)(Z,U) = - S(R(X,Y)Z,U) - S(Z,R(X,Y)U),
\end{equation}
%%%%%%%%%%%%%%%%%
\begin{equation}
\label{eqn1.3}
Q(g,S)(Z,U;X,Y) = - S((X\wedge_g Y)Z,U) - S(Z,(X\wedge_g Y)U),
\end{equation}
%%%%%%%%%%%%%%%%%
and
\begin{equation}
\label{eqn1.4}
(X\wedge_g Y)Z = g(Y,Z)X - g(X,Z)Y
\end{equation}
%%%%%%%%%%%%%%%%%
for all vector fields $X, Y, Z, U$ of $M$, $R$ denotes the curvature tensor of $M$.\\
Then $(M^n,g)$ is Ricci pseudosymmetric if and only if
%%%%%%%%%%%%%%%%%%%%%
\begin{equation}
\label{eqn1.5}
(R(X,Y)\cdot S)(Z,U) = L_S Q(g,S)(Z,U;X,Y)
\end{equation}
%%%%%%%%%%%%%%%%%
holds on $U_S = \{x\in M: S - \frac{r}{n}g \neq 0 \text{ at }x\}$,
for some function $L_S$ on $U_S$. If $R\cdot S = 0$, then $M^n$ is
called Ricci semisymmetric. Every Ricci semisymmetric manifold is
Ricci pseudosymmetric but the converse is not true~\cite{DESZCZ1}.
In this connection it is mentioned that Hui et. al (\cite{SKH-RSL},
\cite{HUI1}) studied Ricci pseudosymmetric generalized
quasi-Einstein manifolds.\\
%%%%%%%%%%%%%%%%%%%
\indent Motivated by the above studies, the object of the present
paper is to study Ricci pseudosymmetric $(LCS)_n$-manifolds whose
metric is a Ricci soliton. In this connection it is mentioned that Hui and Chakraborty \cite{SKH6} 
studied Ricci almost solitons on concircular Ricci pseudosymmetric $\beta$-Kenmotsu manifolds. 
The paper is organized as follows. Section~2 is concerned with preliminaries. In section~3, we
investigate, Ricci solitons on concircular Ricci pseudosymmetric
$(LCS)_n$-manifolds, projective Ricci pseudosymmetric
$(LCS)_n$-manifolds, $W_{3}$-Ricci pseudosymmetric
$(LCS)_n$-manifolds, conharmonic Ricci pseudosymmetric
$(LCS)_n$-manifolds, conformal Ricci pseudosymmetric
$(LCS)_n$-manifolds respectively. Here each curvature tensor has
geometrical significance and hence each type of Ricci
pseudosymmetries has different geometrical interpretance. In each of
the cases, we found the value of $L_S$ and hence it turns out that
the condition that a Ricci soliton is shrinking, steady, or
expanding depends on $L_S$ being less than, equal, or greater than
certain value. We call it the \textsl{Critical Value} for $L_S$. In
each type of Ricci pseudosymmetry, the critical value for $L_S$ is
obtained. Finally we construct an example of concircular Ricci
pseudosymmetric $(LCS)_3$-manifold whose metric is Ricci soliton
through which Theorem 3.1 is verified.
%%%%%%%%%%%%%%%%%%%%%%%%%%%%%%%%%%%%%%%%%%%%%%%%%%%%%%%%%%%%%%%%%%%%%%%%%%%%%%%%%%%%%%%%%%%%%5%%%%%%%%%%%%%%%%%%%%%
\section{Preliminaries}
%%%%%%%%%%%%%%%%%
An $n$-dimensional Lorentzian manifold $M$ is a smooth connected
paracompact Hausdorff manifold with a Lorentzian metric $g$, that
is, $M$ admits a  smooth symmetric tensor field $g$ of type (0, 2)
such that for each point $p\in M$, the tensor $g_{p}:T_{p}M\times
T_{p}M$ $\rightarrow\mathbb{R}$ is a non-degenerate inner product of
signature $(-, +, \cdots, +)$, where $T_{p}M$ denotes the tangent
vector space of $M$ at $p$. A non-zero vector $v$ $\in T_{p} M$ is
said to be timelike (resp. null, spacelike) if it satisfies
$g_{p}(v,v) < 0$ (resp. = 0, $> 0$) \cite{NIL1}.
%%%%%%%%%%%%%%%%%%%%%%%%%%%%
\begin{definition}
(\cite{SHAIKH1},\cite{SHAIKH3}) In a Lorentzian manifold $(M, g)$ a
vector field $P$ defined by
%%%%%%%%%
\begin{equation*}
g(X,P)=A(X),
\end{equation*}
%%%%%%%%%%%%%
for any $X\in\Gamma(TM)$, the section of all smooth tangent vector fields on $M$,
 is said to be a concircular vector field if
%%%%%%%%%5
\begin{equation*}
({\nabla}_{X}A)(Y)=\alpha \{g(X,Y)+\omega(X)A(Y)\}
\end{equation*}
%%%%%%%%%%%%
where $\alpha$ is a non-zero scalar function and $\omega$ is a
closed 1-form and ${\nabla}$ denotes the operator of covariant
differentiation with respect to the Lorentzian metric $g$.
\end{definition}
%%%%%%%%%%%%%%%%%%%%%%%%%%%
Let $M$ be an $n$-dimensional Lorentzian manifold admitting a unit
timelike concircular vector field $\xi$, called the characteristic
vector field of the manifold. Then we have
%%%%%%%%%5
\begin{equation}
\label{2.1} g(\xi, \xi)=-1.
\end{equation}
%%%%%%%%%%%
Since $\xi$ is a unit concircular vector field, it follows that
there exists a non-zero 1-form $\eta$ such that for
%%%%%%%%%5
\begin{equation}
\label{2.2} g(X,\xi)=\eta(X),
\end{equation}
%%%%%%%%%%%%%%
the equation of the following form holds
%%%%%%%%%5
\begin{equation}
\label{2.3} (\nabla_{X}\eta)(Y)=\alpha \{g(X,Y)+\eta(X)\eta(Y)\}, \
\ \ \alpha\neq 0,
\end{equation}
%%%%%%%%%%
and consequently, we get
%%%%%%%%%5
\begin{equation}
\label{2.4} \nabla_{X}\xi = \alpha [X + \eta(X)\xi]
\end{equation}
%%%%%%%%%%
for all vector fields $X$, $Y$, where $\alpha$ satisfying
%%%%%%%%%
\begin{equation}
\label{2.5} {\nabla}_{X}\alpha =(X\alpha) =d\alpha(X)=\rho\eta(X),
\end{equation}
%%%%%%%%%%%%%%
$\rho$ being a certain scalar function given by $\rho=-(\xi\alpha)$.
If we put
%%%%%%%%%5
\begin{equation}
\label{2.6} \phi X=\frac{1}{\alpha}\nabla_{X}\xi,
\end{equation}
%%%%%%%%%%%%%
then from (\ref{2.4}) and (\ref{2.6}) we have
%%%%%%%%%5
\begin{equation}
\label{2.7} \phi X= X+\eta(X)\xi,
\end{equation}
%%%%%%%%%%%
from which it follows that $\phi$ is a symmetric (1, 1) tensor and
called the structure tensor of the manifold. Thus the Lorentzian
manifold $M$ together with the unit timelike concircular vector
field $\xi$, its associated 1-form $\eta$ and an (1, 1) tensor field
$\phi$ is said to be a Lorentzian concircular structure manifold
(briefly, $(LCS)_{n}$-manifold), \cite{SHAIKH11}. In particular, if
we take $\alpha = 1$, then we can obtain the LP-Sasakian structure
of Matsumoto \cite{MAT}. In a $(LCS)_{n}$-manifold $(n > 2)$, the
following relations hold
(\cite{SHAIKH11},\cite{SHAIKH3},\cite{SHAIKH4},\cite{SABINA1}):
%%%%%%%%%5
\begin{equation}
\label{2.8} \eta(\xi) = -1,\phi \xi=0,\eta(\phi X)=0,g(\phi X, \phi
Y) = g(X,Y)+\eta(X)\eta(Y),
\end{equation}
%%%%%%%%%5
\begin{equation}
\label{2.9} \phi^2 X= X+\eta(X)\xi,
\end{equation}
%%%%%%%%%5
\begin{equation}
\label{2.10} S(X,\xi)=(n-1)(\alpha^{2}-\rho)\eta(X),
\end{equation}
%%%%%%%%%5
\begin{equation}
\label{2.11} R(X,Y)\xi=(\alpha^{2}-\rho)[\eta(Y)X-\eta(X)Y],
\end{equation}
%%%%%%%%%5
\begin{equation}
\label{2.12} R(\xi,Y)Z=(\alpha^{2}-\rho)[g(Y,Z)\xi-\eta(Z)Y],
\end{equation}
%%%%%%%%%
\begin{equation}
\label{2.13}
({\nabla}_{X}\phi)Y=\alpha\{g(X,Y)\xi+2\eta(X)\eta(Y)\xi+\eta(Y)X\},
\end{equation}
%%%%%%%%%
\begin{equation}
\label{2.14} (X\rho)=d\rho(X)=\beta\eta(X),
\end{equation}
%%%%%%%%%%%%%%
\begin{equation}
\label{2.15} R(X,Y)Z =\phi R(X,Y)Z + (\alpha^2 -
\rho)\big\{g(Y,Z)\eta(X) - g(X,Z)\eta(Y)\big\}\xi,
\end{equation}
%%%%%%%%%%%%
\begin{equation}
\label{2.16} S(\phi X,\phi Y) = S(X,Y) + (n-1)(\alpha^2 -
\rho)\eta(X)\eta(Y)
\end{equation}
%%%%%%%%%%%%
for any vector fields $X,\ Y,\ Z$ on $M$ and $\beta = -(\xi\rho)$
is a scalar function, where $R$ is the curvature tensor and $S$ is the Ricci tensor
of the manifold.\\
%%%%%%%%%%%%%%%%%%%%%55%%%%%%%%%%%%%%%%%%%%%%%%%%%%%%%%%%%%%%%%%%%%%%%%%5
\indent Let $(g, \xi, \lambda)$ be a Ricci soliton on a
$(LCS)_n$-manifold $M$. From (\ref{2.4}), we get
\begin{eqnarray*}
(\pounds_\xi g)(X,Y)&=& g(\nabla_X \xi, Y) + g(X , \nabla_Y \xi)\\
&=& \alpha [g(X + \eta(X)\xi, Y) + g(X, Y+\eta(Y) \xi)]\notag \\
 &=& 2\alpha[g(X,Y)+\eta(X)\eta(Y)],\notag
\end{eqnarray*}
i.e.
\begin{equation}
\label{eqn2.17} \frac{1}{2}(\pounds_\xi g)(X,Y) ={\alpha}\{g(X,Y) +
\eta(X)\eta(Y)\}.
\end{equation}
%%%%%%%%%%%%
From (\ref{eqn1.1}) and (\ref{eqn2.17}) we have
%%%%%%%%%%%%%%%%%%%
\begin{equation}
\label{eqn2.18} S(X,Y) = -(\alpha + \lambda)g(X,Y) -
\alpha\eta(X)\eta(Y),
\end{equation}
%%%%%%%%%%%%%%%%%%%%%%%%%%%%%%%%%%%%%%%%%%%%%%%%%%%%
which yields
\begin{equation}
\label{eqn2.19} QX = -(\alpha + \lambda)X - \alpha\eta(X)\xi,
\end{equation}
%%%%%%%%%%%%%%%%%%%%%%%%%%%%%%%%%%%%%%%%%%%%%%%%%%%%
\begin{equation}
\label{eqn2.20} S(X,\xi) = -\lambda\eta(X),
\end{equation}
%%%%%%%%%%%%%%%%%%%%%%%%%%%%%%
\begin{equation}
\label{eqn2.21} r = -\lambda n - (n-1)\alpha,
\end{equation}
%%%%%%%%%%%%%%%%%%%%%%%%%%%%%%%%%%%%%%%%%%%%%%%%%%%%
where $Q$ is the Ricci operator, i.e., $g(QX,Y) = S(X,Y)$ for all $X$, $Y$ and $r$ is the scalar curvature of $M$.
%%%%%%%%%%%%%%%%%%%%%%%%%%%%%%%%%%%%%%%%%%%%%%%%%%%%%%%%%%%%%%%%%%%%%%%%%%%
\section{Ricci solitons on Ricci pseudosymmetric $(LCS)_n$-manifolds}
%%%%%%%%%%%%%%%%%%%%%%%%%%%%%%%%%%%%
This section deals with the study of Ricci solitons on concircular
(resp., projective, $W_3$, conharmonic, conformal) Ricci
pseudosymmetric $(LCS)_n$- manifolds. A concircular curvature tensor
is an interesting invariant of a concircular transformation. A
transformation of a $(LCS)_n$-manifold $M$, which transforms every
geodesic circle of $M$ into a geodesic circle, is called a
concircular transformation \cite{YANO}. A concircular transformation
is always a conformal transformation \cite{KUH}. Here geodesic
circle means a curve in $M$ whose first curvature is constant and
whose second curvature is identically zero. Thus the geometry of
concircular transformations, that is, the concircular geometry, is a
generalization of inversive geometry in the sense that the change of
metric is more general than that induced by a circle preserving
diffeomorphism. The interesting invariant of a concircular
transformation is the concircular curvature tensor $\tilde{C}$,
which is defined by \cite{YANO}
%%%%%%%%%%%%%%%%%%%%%%%%
\begin{equation}
\label{eqn3.1}
\tilde{C}(X,Y)Z = R(X,Y)Z-\frac{r}{n(n-1)}\big[g(Y,Z)X-g(X,Z)Y\big],
\end{equation}
%%%%%%%%%%%%%%%%%%%%%%
where $R$ is the curvature tensor and $r$ is the scalar curvature of the manifold.
Also $(LCS)_n$-manifolds with vanishing concircular curvature tensor are of constant curvature.
Thus, the concircular curvature tensor is a measure of the failure of a
$(LCS)_n$-manifold to be of constant curvature.\\
Using (\ref{2.2}), (\ref{2.11}) and (\ref{2.15}), we get
%%%%%%%%%%%%%%%%%%%%%%%%%%%%%%%%%5
\begin{equation}
\label{eqn3.2}
\tilde{C}(X,Y)\xi = \large[(\alpha^2 - \rho) - \frac{r}{n(n-1)}\large][\eta(Y)X - \eta(X)Y],
\end{equation}
%%%%%%%%%%%%%%%%%%%%%%%%
\begin{equation}
\label{eqn3.3}
\eta(\tilde{C}(X,Y)U) = [\frac{r}{n(n-1)} - (\alpha^2 - \rho)][\eta(Y)g(X,U) - \eta(X)g(Y,U)].
\end{equation}
%%%%%%%%%%%%%%%%%%%%%%%%%%%%
\indent A $(LCS)_n$-manifold $(M^n,g)$ is said to be concircular Ricci pseudosymmetric if its
concircular curvature tensor $\tilde{C}$ satisfies
%%%%%%%%%%%%%%%%%%%%%
\begin{equation}
\label{eqn3.4}
(\tilde{C}(X,Y)\cdot S)(Z,U) = L_S Q(g,S)(Z,U;X,Y)
\end{equation}
%%%%%%%%%%%%%%%%%
on $U_S = \{x\in M: S\neq \frac{r}{n}g \text{ at }x\}$, where $L_S$
is some function on $U_S$.\\
\indent Let us take a concircular Ricci pseudosymmetric $(LCS)_n$-manifold whose metric is Ricci soliton.
Then by virtue of (\ref{eqn3.4}) that
%%%%%%%%%%%%%%%%%%%%%
\begin{eqnarray}
\label{eqn3.5}
S(\tilde{C}(X,Y)Z,U) + S(Z,\tilde{C}(X,Y)U) = L_S [g(Y,Z)S(X,U)\\
- g(X,Z)S(Y,U) + g(Y,U)S(X,Z) - g(X,U)S(Y,Z)].\notag
\end{eqnarray}
%%%%%%%%%%%%%%%%%
By virtue of (\ref{eqn2.18}) it follows from (\ref{eqn3.5}) that
%%%%%%%%%%%%%%%%%%%%%
\begin{eqnarray}
\label{eqn3.6}
&&\eta(\tilde{C}(X,Y)Z) \eta(U) + \eta(Z) \eta(\tilde{C}(X,Y)U)\\
&=& L_S[g(Y,Z)\eta(X)\eta(U) - g(X,Z)\eta(Y)\eta(U)\notag\\
&+& g(Y,U)\eta(X)\eta(Z) - g(X,U)\eta(Y)\eta(Z)].\notag
\end{eqnarray}
%%%%%%%%%%%%%%%%%
Setting $Z = \xi$ in (\ref{eqn3.6}) and using (\ref{eqn3.2}) and (\ref{eqn3.3}), we get
%%%%%%%%%%%%%%%%%%%%%
\begin{equation}
\label{eqn3.7}
[L_{S} - \{(\alpha^2 - \rho) - \frac{r}{n(n-1)}\}][\eta(Y)g(X,U) - \eta(X)g(Y,U)] = 0.
\end{equation}
%%%%%%%%%%%%%%%%%
Putting $Y=\xi$ in (\ref{eqn3.7}) and using (\ref{2.8}) and
(\ref{eqn2.21}), we get
%%%%%%%%%%%%%%%%%%%%%
\begin{equation}
\label{eqn3.8}
[L_{S} - \{(\alpha^2 - \rho) + \frac{\lambda}{n-1} + \frac{\alpha}{n}\}][g(X,U) + \eta(X)\eta(U)] = 0
\end{equation}
%%%%%%%%%%%%%%%%%
for all vector fields $X$ and $U$, which follows that
%%%%%%%%%%%%%%%%%%%%%%%%%%%%%%%%%%%%%%%%%%%%%%
\begin{equation}
\label{eqn3.9}
L_{S} = (\alpha^2 - \rho) + \frac{\lambda}{n-1} + \frac{\alpha}{n}.
\end{equation}
%%%%%%%%%%%%%%%%%
This leads to the following:
%%%%%%%%%%%%%%%%%%%%%%%
\begin{theorem}
If $(g,\xi,\lambda)$ is a Ricci soliton on a concircular Ricci pseudosymmetric $(LCS)_n$-manifold, then
$L_{S} = (\alpha^2 - \rho) + \frac{\lambda}{n-1} + \frac{\alpha}{n}$.
\end{theorem}
%%%%%%%%%%%%%%%%%%%%%%%%%%%%%
Again from (\ref{eqn3.9}), we get
%%%%%%%%%%%%%%%%%%%%%%%%%%%%%%%%%%%%%%
\begin{equation}
\label{eqn3.10}
\lambda = (n-1)[L_{S} - \{\frac{\alpha}{n} + (\alpha^2 - \rho)\}].
\end{equation}
%%%%%%%%%%%%%%%%%%%%%%%%%%%%%
Since $n>1$, we have from (\ref{eqn3.10}) that $\lambda <0$, $= 0$ and $> 0$ according as $L_{S}<\frac{\alpha}{n} + (\alpha^2 - \rho)$,
$L_{S} = \frac{\alpha}{n} + (\alpha^2 - \rho)$ and $L_{S} > \frac{\alpha}{n} + (\alpha^2 - \rho)$ respectively. This leads to the following:
%%%%%%%%%%%%%%%%%%%%%%%%%%%%%%%%%%
\begin{corollary}
In a concircular Ricci pseudosymmetric $(LCS)_n$-manifold, the
\textsl{Critical Value} for $L_{S}$ is $\frac{\alpha}{n} + (\alpha^2 -
\rho)$.
\end{corollary}
%%%%%%%%%%%%%%%%%%%%%%%%%%%%%%%%%%%%%%%%%%%%%%%%%%%%%%%
\begin{example}
We consider a $3$-dimensional manifold $M=\{(x,y,z)\in
\mathbb{R}^3:z\neq 0\}$, where $(x,y,z)$ are the standard
coordinates in $\mathbb{R}^3$. Let $\{E_1,E_2,E_3\}$ be a linearly
independent global frame on $M$ given by
$$E_1= z^2\frac{\partial}{\partial x},E_2= z^2\frac{\partial}{\partial y},E_3=\frac{\partial}{\partial z}.$$
Let $g$ be the Lorentzian metric defined by
$g(E_1,E_2)$ = $g(E_1,E_3)$ = $g(E_2,E_3)$ = 0, $g(E_1,E_1) = g(E_2,E_2) = 1, g(E_3,E_3) = -1$.
Let $\eta$ be the $1$-form defined by $\eta(U)=g(U,E_3)$ for any
$U\in\chi(M)$. Let $\phi$ be the $(1,1)$ tensor field defined by
$\phi E_1 = E_1, \phi E_2 = E_2$ and $\phi E_3 = 0$. Then using the
linearity of $\phi$ and $g$ we have
$$\eta(E_3)= -1, \phi U = U + \eta(U)E_3$$ and $$g(\phi U, \phi W) = g(U,W) + \eta(U)\eta(W)$$
for any $U, W \in\chi (M) $. Let $\nabla$ be the Levi-Civita connection
for the Lorentzian metric $g$ and $R$ be the curvature tensor of
$g$. Then we have $$[E_1,E_2]=0, [E_1,E_3]=-\frac{2}{z}E_1,
[E_2,E_3]=-\frac{2}{z}E_2.$$ Using Koszul formula for the Lorentzian
metric $g$, we can easily calculate
$$\nabla_{E_1}E_1 = -\frac{2}{z}E_3, \nabla_{E_1}E_2 = 0, \nabla_{E_1}E_3 = -\frac{2}{z}E_1,$$
$$\nabla_{E_2}E_1 = 0, \nabla_{E_2}E_2 = -\frac{2}{z}E_3, \nabla_{E_2}E_3 = -\frac{2}{z}E_2,$$
$$\nabla_{E_3}E_1 =,\nabla_{E_3}E_2 =,\nabla_{E_3}E_3 = 0.$$
From the above, it can be easily seen that for $E_3 = \xi,
(\phi,\xi,\eta,g)$ is a $(LCS)_3$ structure on $M$. Consequently
$M^3(\phi,\xi,\eta,g)$ is a $(LCS)_3$-manifold with $\alpha =
-\frac{2}{z}\neq 0$, such that $X(\alpha)=\rho \eta(X)$, where $\rho
= -\frac{2}{z^2}$\cite{HUI}. Using the above relations, we can
easily calculate the non-vanishing components of the curvature
tensor as follows:
%%%%%%%%%%%%%%%%%%%%%%%%%%%
\begin{equation*}
R(E_1,E_2)E_1 = -\frac{4}{z^2}E_2, R(E_1,E_2)E_2 = \frac{4}{z^2}E_1,
\end{equation*}
\begin{equation*}
R(E_1,E_3)E_1 = -\frac{6}{z^2}E_3, R(E_1,E_3)E_3 = -\frac{6}{z^2}E_1,
\end{equation*}
\begin{equation*}
R(E_2,E_3)E_2 = -\frac{6}{z^2}E_3, R(E_2,E_3)E_3 = -\frac{6}{z^2}E_2
\end{equation*}
%%%%%%%%%%%%%%%%%%%%%%%%%%%%%%%%%%%%%%%%%5
and the components which can be obtained from these by the symmetry
properties from which, we can easily calculate the non-vanishing
components of the Ricci tensor as follows:
\begin{equation*}
S(E_1,E_1)= S(E_2,E_2)= -\frac{2}{z^2}, S(E_3,E_3)= -\frac{12}{z^2}.
\end{equation*}
%%%%%%%%%%%%%%%%%%%%%%%%%%%%%%%%%%%
Also, the scalar curvature $r$ is given by:
\begin{eqnarray*}
r&=&\sum_{i=1}^{3} g(E_i,E_i)S(E_i,E_i)\\
&=& S(E_1,E_1)+S(E_2,E_2)-S(E_3,E_3)\\
&=& \frac{8}{z^2}.\notag\\
\end{eqnarray*}
Since $\{E_1,E_2,E_3\}$ forms a basis of the $(LCS)_3$-manifold, any
vector field $X,Y,Z,U \in \chi(M)$ can be written as
\begin{equation*}
X=a_1E_1+b_1E_2+c_1E_3,
\end{equation*}
\begin{equation*}
Y=a_2E_1+b_2E_2+c_2E_3,
\end{equation*}
\begin{equation*}
Z=a_3E_1+b_3E_2+c_3E_3,
\end{equation*}
\begin{equation*}
U=a_4E_1+b_4E_2+c_4E_3,
\end{equation*}
where $a_i,b_i,c_i \in \mathbb{R}^+$ for all $i=1,2,3$ such that
$a_i,b_i,c_i$ are not proportional. Then
\begin{eqnarray}
\label{eqn3.11}
R(X,Y)Z&=&\frac{2}{z^2}\{2b_3(a_1b_2-a_2b_1)-3c_3(a_1c_2-a_2c_1)\}E_1\\
&-&\frac{2}{z^2}\{2a_3(a_1b_2-a_2b_1)+3c_3(b_1c_2-b_2c_1)\}E_2\notag\\
&-&\frac{6}{z^2}\{b_3(b_1c_2-b_2c_1)+a_3(a_1c_2-a_2c_1)\}E_3,\notag
\end{eqnarray}
%%%%%%%%%%%%%%%%%%%%%%%%%%%%%%%%%555
\begin{eqnarray}
\label{eqn3.11ab}
R(X,Y)U&=&\frac{2}{z^2}\{2b_4(a_1b_2-a_2b_1)-3c_4(a_1c_2-a_2c_1)\}E_1\\
&-&\frac{2}{z^2}\{2a_4(a_1b_2-a_2b_1)+3c_4(b_1c_2-b_2c_1)\}E_2\notag\\
&-&\frac{6}{z^2}\{b_4(b_1c_2-b_2c_1)+a_4(a_1c_2-a_2c_1)\}E_3.\notag
\end{eqnarray}
In view of (\ref{eqn3.11}) we have from (\ref{eqn3.1}) that
%%%%%%%%%%%%%%%%%%%%%%%%%%%%%%%%%%%%%%%%%%%%%%%%
\begin{eqnarray*}
\tilde{C}(X,Y)Z&=&R(X,Y)Z-\frac{r}{6}[g(Y,Z)X-g(X,Z)Y]\\
&=&\frac{2}{z^2}[2b_3(a_1b_2-a_2b_1)-3c_3(a_1c_2-a_2c_1)\notag\\
&-&\frac{2}{3}\{a_1(b_2b_3-c_2c_3)-a_2(b_1b_3-c_3c_1)\}]E_1\notag\\
&-&\frac{2}{z^2}[2a_3(a_1b_2-a_2b_1)+3c_3(b_1c_2-b_2c_1)\notag\\
&+&\frac{2}{3}\{b_1(a_2a_3-c_2c_3)-b_2(a_1a_3-c_3c_1)\}]E_2\notag\\
&-&\frac{2}{z^2}[3\{b_3(b_1c_2-b_2c_1)+a_3(a_1c_2-a_2c_1)\}\notag\\
&+&\frac{2}{3}\{c_1(a_2a_3+b_2b_3)-c_2(a_1a_3+b_1b_3)\}]E_3.\notag
\end{eqnarray*}
%%%%%%%%%%%%%%%%%%%%%%%%%%%%%%%%%%%%%
and hence
%%%%%%%%%%%%%%%%%%%%%%%%%%%%%%%%%%%%%%%%%%%
\begin{eqnarray}
\label{eqn3.12}
&&S(\tilde{C}(X,Y)Z,U)\\
=&-&\frac{4a_4}{z^4}[2b_3(a_1b_2-a_2b_1)-3c_3(a_1c_2-a_2c_1)\notag\\&-&\frac{2}{3}\{a_1(b_2b_3-c_2c_3)-a_2(b_1b_3-c_3c_1)\}]\notag\\
&+&\frac{4b_4}{z^4}[2a_3(a_1b_2-a_2b_1)+3c_3(b_1c_2-b_2c_1)\notag\\&+&\frac{2}{3}\{b_1(a_2a_3-c_2c_3)-b_2(a_1a_3-c_3c_1)\}]\notag\\
&+&\frac{24c_4}{z^4}[3\{b_3(b_1c_2-b_2c_1)+a_3(a_1c_2-a_2c_1)\}\notag\\&+&\frac{2}{3}\{c_1(a_2a_3+b_2b_3)-c_2(a_1a_3+b_1b_3)\}].\notag
\end{eqnarray}
Similarly we obtain,
\begin{eqnarray}
\label{eqn3.13} &&S(Z,\tilde{C}(X,Y)U)\\
&=&-\frac{4a_3}{z^4}[2b_4(a_1b_2-a_2b_1)-3c_4(a_1c_2-a_2c_1)\notag\\&-&\frac{2}{3}\{a_1(b_2b_4-c_2c_4)-a_2(b_1b_4-c_1c_4)\}]\notag\\
&+&\frac{4b_3}{z^4}[2a_4(a_1b_2-a_2b_1)+3c_4(b_1c_2-b_2c_1)\notag\\&+&\frac{2}{3}{b_1(a_2a_4-c_2c_4)-b_2(a_1a_4-c-1c-4)}]\notag\\
&+&\frac{24c_3}{z^4}[3\{b_4(b_1c_2-b_2c_1)+a_4(a_1c_2-a_2c_1)\}\notag\\&+&\frac{2}{3}\{c_1(a_2a_4+b_2b_4)-c_2(a_1a_4+b_1b_4)\}].\notag
\end{eqnarray}
Now we have
\begin{eqnarray}
\label{eqn3.14}
\begin{cases}
g(Y,Z)=a_2a_3+b_2b_3-c_2c_3,\\
g(X,Z)=a_1a_3+b_1b_3-c_1c_3,\\
g(Y,U)=a_2a_4+b_2b_4-c_2c_4,\\
g(X,U)=a_1a_4+b_1b_4-c_1c_4.
\end{cases}
\end{eqnarray}
Also we have
\begin{eqnarray}
\label{eqn3.15}
\begin{cases}
S(Y,Z)=-\frac{2}{z^2}(a_2a_3+b_2b_3+6c_2c_3),\\
S(X,Z)=-\frac{2}{z^2}(a_1a_3+b_1b_3+6c_1c_3),\\
S(Y,U)=-\frac{2}{z^2}(a_2a_4+b_2b_4+6c_2c_4),\\
S(X,U)=-\frac{2}{z^2}(a_1a_4+b_1b_4+6c_1c_4).
\end{cases}
\end{eqnarray}
Therefore, from (\ref{eqn3.14}) and (\ref{eqn3.15}), we have
%%%%%%%%%%%%%%%%%%%%%%%%%%%%%%%%%%%%%%%%%%%%%
\begin{eqnarray}
\label{eqn3.16}
&&g(Y,Z)S(X,U)-g(X,Z)S(Y,U)\\&+&g(Y,U)S(X,Z)-g(X,U)S(Y,Z)\notag\\
&=&\frac{14}{z^2}[(a_1c_2-a_2c_1)(a_3c_4+a_4c_3)+(b_1c_2-b_2c_1)(b_3c_4+b_4c_3)]\notag\\
&\neq& 0,\notag
\end{eqnarray}
%%%%%%%%%%%%%%%%%%%%%%%%%%%%%%%%%%%
since $a_i, b_i, c_i$ are not proportional and assume that
$(a_1c_2-a_2c_1)(a_3c_4+a_4c_3)+(b_1c_2-b_2c_1)(b_3c_4+b_4c_3) \neq0$.\\
Also from (\ref{eqn3.12}) and (\ref{eqn3.13}) we get
%%%%%%%%%%%%%%%%%%%%%%%%%%%%%%%%%%%%%%%%%%%%
\begin{eqnarray}
\label{eqn3.17}
&&S(\tilde{C}(X,Y)Z,U)+S(Z,\tilde{C}(X,Y)U)\\
&=&\frac{196}{3z^4}[(a_1c_2-a_2c_1)(a_3c_4+a_4c_3)+(b_1c_2-b_2c_1)(b_3c_4+b_4c_3)]\notag\\
&\neq& 0.\notag
\end{eqnarray}
%%%%%%%%%%%%%%%%%%%%%%%%%%%%%%%%%%%%%%%%%%%%%%%%
Let us consider the function
\begin{equation}
\label{eqn3.18}
L_S=\frac{14}{3z^2}.
\end{equation}
%%%%%%%%%%%%%%%%%%%%%%%%%%%%%%%%%%
By virtue of (\ref{eqn3.18}) we have from (\ref{eqn3.16}) and (\ref{eqn3.17}) that
\begin{eqnarray*}
&&S(\tilde{C}(X,Y)Z,U)+S(Z,\tilde{C}(X,Y)U)=L_S[g(Y,Z)S(X,U)\\
&-&g(X,Z)S(Y,U)+g(Y,U)S(X,Z)-g(X,U)S(Y,Z)].\notag
\end{eqnarray*}
%%%%%%%%%%%%%%%%%%%%%%%%%%%%%%%%%%%%%%%%%%%
Hence the $(LCS)_3$-manifold under consideration is concircular
Ricci pseudosymmetric. If $(g,\xi,\lambda)$ is a Ricci soliton on
this $(LCS)_3$-manifold, then from (\ref{eqn2.21}) we get
$$r=-3\lambda-2\alpha,$$
i.e., $$\frac{8}{z^2}=-3\lambda+\frac{4}{z},$$
i.e.,$$\lambda=\frac{4}{3}\big(\frac{1}{z}-\frac{2}{z^2}\big)$$ and hence
from (\ref{eqn3.9}) we get
%%%%%%%%%%%%%%%%%%%%%%%%%%%%%%%%%%%%%%%%%%%%%%%%%
\begin{equation*}
L_S=(\alpha^2-\rho)+\frac{\lambda}{2}+\frac{\alpha}{3}=\frac{14}{3z^2}, \ \ \text{ as } \alpha = -\frac{2}{z}, \rho = -\frac{2}{z^2},
\end{equation*}
%%%%%%%%%%%%%%%%%%%%%%%%%%%%%%%%%5
which satisfies (\ref{eqn3.18}). Thus Theorem 3.1 is verified.
\end{example}
%%%%%%%%%%%%%%%%%%%%%%%%%%%%%%%%%%%%%%%%%%%%%%%%%%%%%%%%%%%%%%%%%%%%%%%%
\indent Now we study of Ricci solitons on projective Ricci
pseudosymmetric $(LCS)_n$-manifolds. The projective curvature tensor
is an important  concept of Riemannian geometry, which one uses to
calculate the basic geometric measurements on a manifold, namely,
angle, distance and various invariants on it. The projective
transformation on a $(LCS)_n$-manifold $(n>1)$ is a transformation
under which geodesic transforms into geodesic. The Weyl projective
curvature tensor is given by \cite{DE-SHAIKH}
%%%%%%%%%%%%%%%%%%%%%%%%
\begin{equation}
\label{eqn4.1}
P(X,Y)Z = R(X,Y)Z-\frac{1}{n-1}\big[S(Y,Z)X-S(X,Z)Y\big],
\end{equation}
%%%%%%%%%%%%%
where $R$ and $S$ are the curvature tensor and Ricci tensor of the
manifold respectively. Using (\ref{2.2}), (\ref{2.11}), (\ref{2.15})
and (\ref{eqn2.18}), we get
%%%%%%%%%%%%%%%%%%%%%%%%%%%%%%%%%5
\begin{equation}
\label{eqn4.2}
P(X,Y)\xi = [(\alpha^2 - \rho) + \frac{\lambda}{n-1}][\eta(Y)X - \eta(X)Y],
\end{equation}
%%%%%%%%%%%%%%%%%%%%%%%%
\begin{equation}
\label{eqn4.3}
\eta(P(X,Y)U) = [(\alpha^2 - \rho) - \frac{\alpha + \lambda}{n-1}][\eta(Y)g(X,U) - \eta(X)g(Y,U)].
\end{equation}
%%%%%%%%%%%%%%%%%%%%%%%%%%%%
\indent A $(LCS)_n$-manifold $(M^n,g)$ is said to be projective Ricci pseudosymmetric if its
projective curvature tensor $P$ satisfies
%%%%%%%%%%%%%%%%%%%%%
\begin{equation}
\label{eqn4.4}
(P(X,Y)\cdot S)(Z,U) = L_S Q(g,S)(Z,U;X,Y).
\end{equation}
%%%%%%%%%%%%%%%%%
holds on $U_S = \{x\in M: S\neq \frac{r}{n}g \text{ at }x\}$, where $L_S$
is some function on $U_S$.\\
\indent Let us take a projective Ricci pseudosymmetric $(LCS)_n$-manifold whose metric is Ricci soliton.
Then we get from (\ref{eqn4.4}) that
%%%%%%%%%%%%%%%%%%%%%
\begin{eqnarray}
\label{eqn4.5}
S(P(X,Y)Z,U) + S(Z,P(X,Y)U) = L_S [g(Y,Z)S(X,U)\\
- g(X,Z)S(Y,U) + g(Y,U)S(X,Z) - g(X,U)S(Y,Z)].\notag
\end{eqnarray}
%%%%%%%%%%%%%%%%%
Using (\ref{eqn2.18}) in (\ref{eqn4.5}), we get
%%%%%%%%%%%%%%%%%%%%%
\begin{eqnarray}
\label{eqn4.6}
&&(\alpha+\lambda)[g(P(X,Y)Z,U) + g(Z,P(X,Y)U)]\\
&+&\alpha[\eta(P(X,Y)Z) \eta(U) + \eta(Z) \eta(P(X,Y)U)]\notag\\
&=& \alpha L_S[g(Y,Z)\eta(X)\eta(U) - g(X,Z)\eta(Y)\eta(U)\notag\\
&+& g(Y,U)\eta(X)\eta(Z) - g(X,U)\eta(Y)\eta(Z)].\notag
\end{eqnarray}
%%%%%%%%%%%%%%%%%
Setting $Z = \xi$ in (\ref{eqn4.6}), we get
%%%%%%%%%%%%%%%%%%%%%
\begin{equation}
\label{eqn4.7}
[\alpha L_{S} - (\alpha + 2\lambda)(\alpha^2 - \rho)][\eta(Y)g(X,U) - \eta(X)g(Y,U)] = 0.
\end{equation}
%%%%%%%%%%%%%%%%%
Putting $Y=\xi$ in (\ref{eqn4.7}) and using (\ref{2.8}), we get
%%%%%%%%%%%%%%%%%%%%%aa
\begin{equation}
\label{eqn4.8}
[\alpha L_{S} - (\alpha + 2\lambda)(\alpha^2 - \rho)][g(X,U) + \eta(X)\eta(U)] = 0
\end{equation}
%%%%%%%%%%%%%%%%%
for all vector fields $X$ and $U$, from which it follows that
%%%%%%%%%%%%%%%%%%%%%%%%%%%%%%%%%%%%%%%%%%%%%%
\begin{equation}
\label{eqn4.9}
L_{S} = \large(1+\frac{2\lambda}{\alpha}\large)(\alpha^2 - \rho).
\end{equation}
%%%%%%%%%%%%%%%%%
This leads to the following:
%%%%%%%%%%%%%%%%%%%%%%%
\begin{theorem}
If $(g,\xi,\lambda)$ is a Ricci soliton on a projective Ricci pseudosymmetric $(LCS)_n$-manifold, then
$L_{S} = (1+\frac{2\lambda}{\alpha})(\alpha^2 - \rho)$.
\end{theorem}
%%%%%%%%%%%%%%%%%%%%%%%%%%%%%
Again from (\ref{eqn4.9}), we get
%%%%%%%%%%%%%%%%%%%%%%%%%%%%%%%%%%%%%%
\begin{equation}
\label{eqn4.10}
\lambda = \frac{\alpha[L_{S} - (\alpha^2 - \rho)]}{2(\alpha^2 - \rho)}.
\end{equation}
%%%%%%%%%%%%%%%%%%%%%%%%%%%%%
This leads to the following:
%%%%%%%%%%%%%%%%%%%%%%%%%%%%%%%%%%
\begin{corollary}
In a projective Ricci pseudosymmetric $(LCS)_n$-manifold, the
\textsl{Critical Value} for $L_S$ is $(\alpha^2-\rho)$, provided
$\frac{\alpha}{(\alpha^2-\rho)}>0$.
\end{corollary}
%%%%%%%%%%%%%%%%%%%%%%%%%%%%%
\noindent{\bf Remark:} In \cite{ABI} Ashoka, Bagewadi and Ingalahalli studied Ricci
solitons in $LCS)_n$-manifolds satisfying $R(\xi,X). \tilde{P}=0$,
where $\tilde{P}$ is the pseudo projective curvature tensor. Thus
the present result in our paper are not just special cases of
results in \cite{ABI}.\\
%%%%%%%%%%%%%%%%%%%%%%%%%%%%%%%%%%%%%%%%%%%%%%%%%%%%%%%%%%%%%%%%%%%%%%%%%%%%%%%%%
\indent Now we study of Ricci solitons on $W_3$-Ricci
pseudosymmetric $(LCS)_n$-manifolds. In 1973 Pokhariyal \cite{POKH}
introduced the notion of a new curvature tensor, denoted by $W_{3}$
and studied its relativistic significance. The $W_{3}$-curvature
tensor of type (1,3) on a $(LCS)_n$-manifold is defined by
%%%%%%%%%%%%%%%%%%%%%%%%%%%%
\begin{equation}
\label{eqn5.1}
W_{3}(X,Y)Z = R(X,Y)Z + \frac{1}{n-1}\big[g(Y,Z)QX - S(X,Z)Y\big],
\end{equation}
%%%%%%%%%%%%%%%%%%%%%%%%%%%%%
where $R$ is the curvature tensor and $Q$ is the Ricci-operator, i.e., $g(QX,Y) = S(X,Y)$ for all $X$, $Y$.\\
Using (\ref{2.11}), (\ref{2.15}), (\ref{eqn2.18}) and
(\ref{eqn2.19}), we get
%%%%%%%%%%%%%%%%%%%%%%%%%%%%%%%%%5
\begin{eqnarray}
\label{eqn5.2}
W_{3}(X,Y)\xi &=& [(\alpha^2 - \rho) - \frac{\lambda}{n-1}][\eta(Y)X - \eta(X)Y]\\
\nonumber&-&\frac{\alpha}{n-1}\eta(Y)[X+\eta(X)\xi],
\end{eqnarray}
%%%%%%%%%%%%%%%%%%%%%%%%
\begin{eqnarray}
\label{eqn5.3}
\ \ \ \ \ \ \ \ \eta(W_{3}(X,Y)U) &=& [(\alpha^2 - \rho) + \frac{\lambda}{n-1}][\eta(Y)g(X,U) - \eta(X)g(Y,U)]\\
\nonumber&+&\frac{\alpha}{n-1}\eta(Y)\{g(X,U)+\eta(X)\eta(U)\}.
\end{eqnarray}
%%%%%%%%%%%%%%%%%%%%%%%%%%%%
\indent A $(LCS)_n$-manifold $(M^n,g)$ is said to be $W_3$-Ricci pseudosymmetric if it satisfies
%%%%%%%%%%%%%%%%%%%%%%%%%%
\begin{equation}
\label{eqn5.4}
(W_{3}(X,Y)\cdot S)(Z,U) = L_S Q(g,S)(Z,U;X,Y)
\end{equation}
%%%%%%%%%%%%%%%%%
holds on $U_S = \{x\in M: S\neq \frac{r}{n}g \text{ at }x\}$, where $L_S$
is some function on $U_S$.\\
\indent Let us take a $W_{3}$-Ricci pseudosymmetric $(LCS)_n$-manifold whose metric is Ricci soliton.
Then we have from (\ref{eqn5.4}) that
%%%%%%%%%%%%%%%%%%%%%
\begin{eqnarray}
\label{eqn5.5}
&&S(W_{3}(X,Y)Z,U) + S(Z,W_{3}(X,Y)U) = L_S [g(Y,Z)S(X,U)\\
&&- g(X,Z)S(Y,U) + g(Y,U)S(X,Z) - g(X,U)S(Y,Z)].\notag
\end{eqnarray}
%%%%%%%%%%%%%%%%%
Using (\ref{eqn2.18}) in (\ref{eqn5.5}), we get
%%%%%%%%%%%%%%%%%%%%%
\begin{eqnarray}
\label{eqn5.6}
&&(\alpha+\lambda)[g(W_{3}(X,Y)Z,U) + g(Z,W_{3}(X,Y)U)]\\
&+&\alpha[\eta(W_{3}(X,Y)Z) \eta(U) + \eta(Z) \eta(W_{3}(X,Y)U)]\notag\\
&=& \alpha L_S[g(Y,Z)\eta(X)\eta(U) - g(X,Z)\eta(Y)\eta(U)\notag\\
&+& g(Y,U)\eta(X)\eta(Z) - g(X,U)\eta(Y)\eta(Z)].\notag
\end{eqnarray}
%%%%%%%%%%%%%%%%%
Setting $Z = \xi$ in (\ref{eqn5.6}) and using (\ref{eqn5.2}) and (\ref{eqn5.3}), we get
%%%%%%%%%%%%%%%%%%%%%
\begin{eqnarray}
\label{eqn5.7}
&&[(\alpha + 2\lambda)(\alpha^2 - \rho) - \frac{\alpha\lambda}{n-1} - \alpha L_{S}][\eta(Y)g(X,U) - \eta(X)g(Y,U)]\\
\nonumber&&- \frac{\alpha^2}{n-1}\eta(Y)[g(X,U) + \eta(X)\eta(U)] = 0.
\end{eqnarray}
%%%%%%%%%%%%%%%%%
Putting $Y=\xi$ in (\ref{eqn5.7}) and using (\ref{2.8}), we get
%%%%%%%%%%%%%%%%%%%%%aa
\begin{equation}
\label{eqn5.8}
[L_{S} - (1 + \frac{2\lambda}{\alpha})(\alpha^2 - \rho) + \frac{\alpha+\lambda}{n-1}][g(X,U) + \eta(X)\eta(U)] = 0
\end{equation}
%%%%%%%%%%%%%%%%%
for all vector fields $X$ and $U$, which follows that
%%%%%%%%%%%%%%%%%%%%%%%%%%%%%%%%%%%%%%%%%%%%%%
\begin{equation}
\label{eqn5.9}
L_{S} = \large(1 + \frac{2\lambda}{\alpha}\large)(\alpha^2 - \rho) - \frac{\alpha+\lambda}{n-1}.
\end{equation}
%%%%%%%%%%%%%%%%%
This leads to the following:
%%%%%%%%%%%%%%%%%%%%%%%
\begin{theorem}
If $(g,\xi,\lambda)$ is a Ricci soliton on a $W_3$-Ricci pseudosymmetric $(LCS)_n$-manifold, then
$L_{S}$ is given by \emph{(\ref{eqn5.9})}.
\end{theorem}
%%%%%%%%%%%%%%%%%%%%%%%%%%%%%
Again from (\ref{eqn5.9}), we get
%%%%%%%%%%%%%%%%%%%%%%%%%%%%%%%%%%%%%%
\begin{equation}
\label{eqn5.10}
\lambda = \frac{(n-1)\alpha}{2(n-1)(\alpha^2 - \rho) - \alpha}[L_{S} - \frac{(n-1)(\alpha^2 - \rho) - \alpha}{n-1}].
\end{equation}
%%%%%%%%%%%%%%%%%%%%%%%%%%%%%
This leads to the following:
%%%%%%%%%%%%%%%%%%%%%%%%%%%%%%%%%%
\begin{corollary}
In a $W_3$-Ricci pseudosymmetric $(LCS)_n$-manifold, the critical
value for $L_S$ is $(\alpha^2-\rho-\frac{\alpha}{n-1})$, provided
$\frac{(n-1)\alpha}{2(n-1)(\alpha^2-\rho)-\alpha}>0$.
\end{corollary}
%%%%%%%%%%%%%%%%%%%%%%%%%%%%%%%%%%%%%%%%%%%%%%%%%%%%%%%%%%%%%%%%%%%%%%%%
\indent We now study of Ricci solitons on conharmonic Ricci
pseudosymmetric $(LCS)_n$-manifolds. Of considerable interest is a
special type of conformal transformations, conharmonic
transformations, which are preserving the harmonicity property of
smooth functions. This type of transformation was introduced by
Ishii \cite{ISHI} in 1957 and is now studied from various points of
view. It is well known that such transformations have a tensor
invariant, the so-called conharmonic curvature tensor. It is easy to
verify that this tensor is an algebraic curvature tensor; that is,
it possesses the classical symmetry properties of the Riemannian
curvature tensor. It is known that a harmonic function is defined as
a function whose Laplacian vanishes. A harmonic function is not
invariant, in general. The conditions under which a harmonic
function remains invariant have been studied by Ishii \cite{ISHI}
who introduced the conharmonic transformation as a subgroup of the
conformal transformation. A manifold whose conharmonic curvature
tensor vanishes at every point of the manifold is called
conharmonically flat manifold. Thus this tensor represents the
deviation of the manifold from canharmonic flatness. As a special
subgroup of the conformal transformation group, Ishii~\cite{ISHI}
introduced the notion of conharmonic transformation under which a
harmonic function transform into a harmonic function. The
conharmonic curvature tensor of type (1,3) on a Riemannian manifold
$(M^n,g)$, $n>3$, is given by \cite{ISHI}.
%%%%%%%%%%%%%%%%%%%%%%%%%%%%%%%%%%%%%%%%%%
\begin{eqnarray}
\label{eqn6.1}
\overline C(X,Y)Z = R(X,Y)Z - \frac{1}{n-2}\big[S(Y,Z)X\\
\nonumber - S(X,Z)Y + g(Y,Z)QX - g(X,Z)QY\big],
\end{eqnarray}
%%%%%%%%%%%%
which is invariant under conharmonic transformation, where $S$ is the Ricci
tensor of the manifold of type (0,2).\\
%%%%%%%%%%%%%%%%%%%%%%%%%%
Using (\ref{2.2}), (\ref{2.11}), (\ref{2.15}) (\ref{eqn2.18}) and
(\ref{eqn2.19}), we get
%%%%%%%%%%%%%%%%%%%%%%%%%%%%%%%%%5
\begin{equation}
\label{eqn6.2}
\overline{C}(X,Y)\xi = [(\alpha^2 - \rho) + \frac{\alpha + 2\lambda}{n-2}][\eta(Y)X - \eta(X)Y],
\end{equation}
%%%%%%%%%%%%%%%%%%%%%%%%
\begin{equation}
\label{eqn6.3}
\eta(\overline{C}(X,Y)U) = [(\alpha^2 - \rho) - \frac{\alpha+2\lambda}{n-2}][\eta(Y)g(X,U) - \eta(X)g(Y,U)].
\end{equation}
%%%%%%%%%%%%%%%%%%%%%%%%%%%%
\indent A $(LCS)_n$-manifold $(M^n,g)$ is said to be conharmonic Ricci pseudosymmetric if its
conharmonic curvature tensor $\overline {C}$ satisfies
%%%%%%%%%%%%%%%%%%%%%
\begin{equation}
\label{eqn6.4}
(\overline {C}(X,Y)\cdot S)(Z,U) = L_S Q(g,S)(Z,U;X,Y).
\end{equation}
%%%%%%%%%%%%%%%%%
holds on $U_S = \{x\in M: S\neq \frac{r}{n}g \text{ at }x\}$, where $L_S$
is some function on $U_S$.\\
\indent Let us take a conharmonic Ricci pseudosymmetric $(LCS)_n$-manifold whose metric is Ricci soliton. Then we get
by virtue of (\ref{eqn6.4}) that
%%%%%%%%%%%%%%%%%%%%%
\begin{eqnarray}
\label{eqn6.5}
S(\overline{C}(X,Y)Z,U) + S(Z,\overline{C}(X,Y)U) = L_S [g(Y,Z)S(X,U)\\
- g(X,Z)S(Y,U) + g(Y,U)S(X,Z) - g(X,U)S(Y,Z)].\notag
\end{eqnarray}
%%%%%%%%%%%%%%%%%
By virtue of (\ref{eqn2.18}) it follows from (\ref{eqn6.5}) that
%%%%%%%%%%%%%%%%%%%%%
\begin{eqnarray}
\label{eqn6.6}
&&\eta(\overline{C}(X,Y)Z) \eta(U) + \eta(Z) \eta(\overline{C}(X,Y)U)\\
&=& L_S[g(Y,Z)\eta(X)\eta(U) - g(X,Z)\eta(Y)\eta(U)\notag\\
&+& g(Y,U)\eta(X)\eta(Z) - g(X,U)\eta(Y)\eta(Z)].\notag
\end{eqnarray}
%%%%%%%%%%%%%%%%%
Setting $Z = \xi$ in (\ref{eqn6.6}) and using (\ref{eqn6.2}) and (\ref{eqn6.3}), we get
%%%%%%%%%%%%%%%%%%%%%
\begin{equation}
\label{eqn6.7}
[L_{S} + (\alpha^2 - \rho) - \frac{\alpha+2\lambda}{n-2}][\eta(Y)g(X,U) - \eta(X)g(Y,U)] = 0.
\end{equation}
%%%%%%%%%%%%%%%%%
Putting $Y=\xi$ in (\ref{eqn6.7}) and using (\ref{2.8}), we get
%%%%%%%%%%%%%%%%%%%%%
\begin{equation}
\label{eqn6.8}
[L_{S} + (\alpha^2 - \rho) - \frac{\alpha+2\lambda}{n-2}][g(X,U) + \eta(X)\eta(U)] = 0
\end{equation}
%%%%%%%%%%%%%%%%%
for all vector fields $X$ and $U$, which follows that
%%%%%%%%%%%%%%%%%%%%%%%%%%%%%%%%%%%%%%%%%%%%%%
\begin{equation}
\label{eqn6.9}
L_{S} = \frac{\alpha+2\lambda}{n-2} - (\alpha^2 - \rho).
\end{equation}
%%%%%%%%%%%%%%%%%
This leads to the following:
%%%%%%%%%%%%%%%%%%%%%%%
\begin{theorem}
If $(g,\xi,\lambda)$ is a Ricci soliton on a conharmonic Ricci pseudosymmetric $(LCS)_n$-manifold, then
$L_{S} = \frac{\alpha+2\lambda}{n-2} - (\alpha^2 - \rho)$.
\end{theorem}
%%%%%%%%%%%%%%%%%%%%%%%%%%%%%
Again from (\ref{eqn6.9}), we get
%%%%%%%%%%%%%%%%%%%%%%%%%%%%%%%%%%%%%%
\begin{equation}
\label{eqn6.10}
\lambda = \frac{1}{2}[(n-2)\{L_{S} + (\alpha^2 - \rho)\} - \alpha].
\end{equation}
%%%%%%%%%%%%%%%%%%%%%%%%%%%%%
Since $n>2$, we have from (\ref{eqn6.10}) that $\lambda <0$, $= 0$
and $> 0$ according as $L_S<, = and >
\frac{\alpha}{n-2}-(\alpha^2-\rho)$, respectively. This leads to the
following:
%%%%%%%%%%%%%%%%%%%%%%%%%%%%%%%%%%
\begin{corollary}
In a conharmonic Ricci pseudosymmetric $(LCS)_n$-manifold,
\textsl{Critical Value} for $L_S$ is
$\frac{\alpha}{n-2}-(\alpha^2-\rho)$.
\end{corollary}
%%%%%%%%%%%%%%%%%%%%%%%%%%%%%%%%%%%%%%%%%%%%%%%%%%%%%%%%%%%%%%%%%%%%%%%%%%%%%%%%%%%%%%%%%%%%%%%%%%%%%%%%%%%%%
\indent Now we study of Ricci solitons on conformal Ricci
pseudosymmetric $(LCS)_n$ manifolds. In differential geometry, the
Weyl curvature tensor, named after Hermann Weyl, is a measure of the
curvature of spacetime or, more generally, a pseudo-Riemannian
manifold. Like the Riemann curvature tensor, the Weyl tensor
expresses the tidal force that a body feels when moving along a
geodesic. The Weyl tensor is the traceless component of the Riemann
curvature tensor \cite{KUH}. Since the trace component of the
Riemann curvature tensor, i.e. the Ricci curvature, contains
precisely the information about how volumes change in the presence
of tidal forces, the Weyl tensor does not convey information on how
the volume of the manifold changes, but rather only how the shape of
the body is distorted by the tidal force.\\
%%%%%%%%%%%%%%%%%%%%%%%%%%%%%%%%%%%%%%%%%%%%%%%%%%%%%%%%%%%%%%%%%%%%%%%%%%%%%%%%%%%%%%%%%%%%%%%%%%%%%%%%%%%55
\indent In general relativity, the Weyl curvature is the only part of the curvature that exists in free space-a solution of the vacuum Einstein
equation-and it governs the propagation of gravitational radiation through regions of space devoid of matter. More generally, the Weyl
curvature is the only component of curvature for Ricci-flat manifolds and always governs the characteristics of the field equations
of an Einstein manifold. In dimensions $2$ and $3$ the Weyl curvature tensor vanishes identically. In dimensions $\geq 4$, the Weyl
curvature is generally nonzero. If the Weyl tensor vanishes in dimension $\geq 4$, then the metric is locally conformally flat: there
exists a local coordinate system in which the metric tensor is proportional to a constant tensor. This fact was a key component
of Nordstr$\ddot{\mbox{o}}$m's theory of gravitation, which was a precursor of general relativity.\\
%%%%%%%%%%%%%%%%%%%%%%%%%%%%%%%%%%%%%%%%%%%%%%%%%%%%%%%%%%5555
\indent The Weyl tensor has the special property that it is
invariant under conformal changes to the metric. For this reason the
Weyl tensor is also called the conformal tensor. It follows that the
necessary condition for a Riemannian manifold to be conformally flat
is that the Weyl tensor vanish. In dimensions $\geq 4$ this
condition is sufficient as well. In dimension $3$ the vanishing of
the Cotton tensor is the necessary and sufficient condition for the
Riemannian manifold being conformally flat. Any $2$-dimensional
(smooth) Riemannian manifold is conformally flat, a consequence of
the existence of isothermal coordinates. Conformal transformations
of a Riemannian structures
are an important object of study in differential geometry.\\
%%%%%%%%%%%%%%%%%%%%%%%%%%%%%%%%%%%%%%%%%%%%%%5555555
\indent The conformal transformation on a $(LCS)_n$-manifold is a transformation
under which the angle between two curves remains invariant. The Weyl conformal
curvature tensor $C$ of type (1,3) of an $n$-dimensional Riemannian manifold
$(LCS)_n$ $(n>3)$ is defined by~\cite{DE-SHAIKH}
%%%%%%%%%%%%%%%%%%%%%
\begin{eqnarray}
\label{eqn7.1}
C(X,Y)Z &=& R(X,Y)Z - \frac{1}{n-2}[S(Y,Z)X - S(X,Z)Y\\
&+&g(Y,Z) QX - g(X,Z)QY]  \notag \\
&+&\frac{r}{(n-1)(n-2)}\{g(Y,Z)X - g(X,Z)Y\},\notag
\end{eqnarray}
%%%%%%%%%%%%%%%%%
where $R$, $S$, $Q$ and $r$ are the Curvature tensor, Ricci tensor,
Ricci-operator and scalar curvature of the manifold respectively.
Using (\ref{2.2}), (\ref{2.11}), (\ref{2.15}) (\ref{eqn2.18}) and
(\ref{eqn2.19}), we get
%%%%%%%%%%%%%%%%%%%%%%%%%%%%%%%%%5
\begin{equation}
\label{eqn7.2}
C(X,Y)\xi = [(\alpha^2 - \rho) + \frac{\lambda}{n-1}][\eta(Y)X - \eta(X)Y],
\end{equation}
%%%%%%%%%%%%%%%%%%%%%%%%
\begin{equation}
\label{eqn7.3}
\eta(C(X,Y)U) = [(\alpha^2 - \rho) - \frac{\lambda}{n-1}][\eta(Y)g(X,U) - \eta(X)g(Y,U)].
\end{equation}
%%%%%%%%%%%%%%%%%%%%%%%%%%%%
\indent A $(LCS)_n$-manifold $(M^n,g)$ is said to be conformal Ricci pseudosymmetric if its
conformal curvature tensor $C$ satisfies
%%%%%%%%%%%%%%%%%%%%%
\begin{equation}
\label{eqn7.4}
(C(X,Y)\cdot S)(Z,U) = L_S Q(g,S)(Z,U;X,Y).
\end{equation}
%%%%%%%%%%%%%%%%%
holds on $U_S = \{x\in M: S\neq \frac{r}{n}g \text{ at }x\}$, where $L_S$
is some function on $U_S$.\\
\indent Let us take a conformal Ricci pseudosymmetric $(LCS)_n$-manifold whose metric is Ricci soliton. Then we have from (\ref{eqn7.4}) that
%%%%%%%%%%%%%%%%%%%%%
\begin{eqnarray}
\label{eqn7.5}
S(C(X,Y)Z,U) + S(Z,C(X,Y)U) = L_S [g(Y,Z)S(X,U)\\
- g(X,Z)S(Y,U) + g(Y,U)S(X,Z) - g(X,U)S(Y,Z)].\notag
\end{eqnarray}
%%%%%%%%%%%%%%%%%
By virtue of (\ref{eqn2.18}) it follows from (\ref{eqn7.5}) that
%%%%%%%%%%%%%%%%%%%%%
\begin{eqnarray}
\label{eqn7.6}
&&\eta(C(X,Y)Z) \eta(U) + \eta(Z) \eta(C(X,Y)U)\\
&=& L_S[g(Y,Z)\eta(X)\eta(U) - g(X,Z)\eta(Y)\eta(U)\notag\\
&+& g(Y,U)\eta(X)\eta(Z) - g(X,U)\eta(Y)\eta(Z)].\notag
\end{eqnarray}
%%%%%%%%%%%%%%%%%
Setting $Z = \xi$ in (\ref{eqn7.6}) and using (\ref{eqn7.2}) and (\ref{eqn7.3}), we get
%%%%%%%%%%%%%%%%%%%%%
\begin{equation}
\label{eqn7.7}
[L_{S} + (\alpha^2 - \rho) - \frac{\lambda}{n-1}][\eta(Y)g(X,U) - \eta(X)g(Y,U)] = 0.
\end{equation}
%%%%%%%%%%%%%%%%%
Putting $Y=\xi$ in (\ref{eqn7.7}) and using (\ref{2.8}), we get
%%%%%%%%%%%%%%%%%%%%%
\begin{equation}
\label{eqn7.8}
[L_{S} + (\alpha^2 - \rho) - \frac{\lambda}{n-1}][g(X,U) + \eta(X)\eta(U)] = 0
\end{equation}
%%%%%%%%%%%%%%%%%
for all vector fields $X$ and $U$, which follows that
%%%%%%%%%%%%%%%%%%%%%%%%%%%%%%%%%%%%%%%%%%%%%%
\begin{equation}
\label{eqn7.9}
L_{S} = \frac{\lambda}{n-1} - (\alpha^2 - \rho).
\end{equation}
%%%%%%%%%%%%%%%%%
This leads to the following:
%%%%%%%%%%%%%%%%%%%%%%%
\begin{theorem}
If $(g,\xi,\lambda)$ is a Ricci soliton on a conformal Ricci pseudosymmetric $(LCS)_n$-manifold, then
$L_{S}$ is given by \emph{(\ref{eqn7.9})}.
\end{theorem}
%%%%%%%%%%%%%%%%%%%%%%%%%%%%%
Again from (\ref{eqn7.9}), we get
%%%%%%%%%%%%%%%%%%%%%%%%%%%%%%%%%%%%%%
\begin{equation}
\label{eqn7.10}
\lambda = (n-1)[L_{S} + (\alpha^2 - \rho)].
\end{equation}
%%%%%%%%%%%%%%%%%%%%%%%%%%%%%
Since $n>1$, we have from (\ref{eqn7.10}) that $\lambda <0$, $= 0$
and $> 0$ according as $L_{S} < -(\alpha^2 - \rho) <0$, $L_{S} =
-(\alpha^2 - \rho) = 0$ and $L_{S} > -(\alpha^2 - \rho) >0$
respectively. This leads to the following:
%%%%%%%%%%%%%%%%%%%%%%%%%%%%%%%%%%
\begin{corollary}
In a conformal Ricci pseudosymmetric $(LCS)_n$-manifold,
 the \textsl{Critical Value} for  $L_{S}$ is $-(\alpha^2 - \rho) $.
\end{corollary}
%%%%%%%%%%%%%%%%%%%%%%%%%%%%%%%%
\section{Summary}
From Theorem 3.1 to Theorem 3.5, we have the following:

\medskip
Let $(g,\xi,\lambda)$ be a Ricci soliton on a $(LCS)_n$ manifold $M$. Then the following holds:
\begin{center}
\begin{tabular}{|c|c|}
\hline $M$ & $L_S$ \\
\hline Concircular Ricci Pseudosymmetric &
$(\alpha^2-\rho)+\frac{\lambda}{n-1}+\frac{\alpha}{n}$\\
\hline Projective Ricci Pseudosymmetric &
$(1+\frac{2\lambda}{\alpha})(\alpha^2-\rho)$\\
\hline $W_3$-Ricci Pseudosymmetric & $(1+\frac{2\lambda}{\alpha})(\alpha^2-\rho)-\frac{\alpha+\lambda}{n-1}$\\
\hline Conharmonic Ricci Pseudosymmetric &
$\frac{\alpha+2\lambda}{n-2}-(\alpha^2-\rho)$\\
\hline Conformal Ricci Pseudosymmetric &
$\frac{\lambda}{n-1}-(\alpha^2-\rho)$\\
\hline
\end{tabular}
\end{center}

Again, from Corollary 3.1 to Corollary 3.5, we have the following:
%%%%%%%%%%%%%%%%%%%%%%%%%%%%%%%%%

\medskip
In a $(LCS)_n$-manifold $M$, the following holds:
\begin{center}
\begin{tabular}{|c|c|}
\hline $M$ & Critical value for $L_S$ \\
\hline Concircular Ricci Pseudosymmetric &
$\frac{\alpha}{n}+(\alpha^2-\rho)$\\
\hline Projective Ricci Pseudosymmetric &
$(\alpha^2-\rho)$, provided $\alpha > \alpha^2-\rho$\\
\hline $W_3$-Ricci Pseudosymmetric & $(\alpha^2-\rho)-\frac{\alpha}{n-1}$, provided $\frac{(n-1)\alpha}{2(n-1)(\alpha^2)-\alpha}>0$\\
\hline Conharmonic Ricci Pseudosymmetric &
$\frac{\alpha}{n-2}-(\alpha^2-\rho)$\\
\hline Conformal Ricci Pseudosymmetric &
$-(\alpha^2-\rho)$\\
\hline
\end{tabular}
\end{center}

\bigskip
%%%%%%%%%%%%%%%%%%%%%%%%%%%%%%%%%%%%%%%%%%
%%%%%%%%%%%%%%%%%%%%%%%%%%%%%%%%%%%%%%%%%%%%%%

%%%%%%%%%%%%%%%%%%%%%%%%%%%%%%%%%%%%%%%%%%%%%%%%%%%%%%%%%%%%%%%%%%%%%%%%%%%%%%
\vspace{0.1in} \noindent S. K. Hui\\
Department of Mathematics\\
The University of Burdwan\\
Burdwan - 713 104\\
West Bengal, India\\
E-mail: shyamal\_hui@yahoo.co.in\\
skhui@math.buruniv.ac.in\\

\noindent R. S. Lemence\\
Institute of Mathematics, College of Science\\
University of Philippines\\
Diliman, Quezon City 1101 Philippines\\
E-mail: rslemence@math.upd.edu.ph\\

\noindent D. Chakraborty\\
Department of Mathematics\\
Sidho Kanho Birsha University\\
Purulia - 723 104\\
West Bengal, India\\
E-mail: debabratamath@gmail.com\\
\end{document}